\documentclass{article}%
\usepackage{amsmath}%
\usepackage{amsfonts}%
\usepackage{amssymb}%
\usepackage{graphicx}

\begin{document}

\title{On Studying the Phase Behavior of the Riemann Zeta Function Along the Critical Line}
\author{Henrik Stenlund}
\date{June 1st, 2018}
\maketitle

\begin{abstract}
The critical line of the Riemann zeta function is studied from a new viewpoint. It is found that the ratio between the zeta function at any zero and the corresponding one at a conjugate point has a certain phase and its absolute value is unity. This fact is valid along the whole critical line and only there. The common  functional equation is used with the aid of the function ratio between any zero and its negative side pair, a complex conjugate. As a result, an equation is obtained for solving the phase along the critical line. 
\footnote{Visilab Report \#2018-06}
\footnote{The author is with Visilab Signal Technologies Oy, Finland. The author is grateful for full support from the organization}
\subsection{Keywords}
Riemann zeta function, zeros of zeta function, recursion relation of zeta function, functional equation of zeta function
\subsection{Mathematical Classification}
MSC: 11M26
\end{abstract}

\section{Introduction}
\subsection{Motivation}
Developing methods for studying of the nature of the nontrivial zeros on the critical line of the Riemann zeta function are very important. One is attempting in this paper to find a way of eliminating the zeta function out of the functional equations allowing possibly new expressions to be developed in terms of more elementary functions. We are attempting to find an expression for the phase of the zeta function along the critical line.
\subsection{Preliminaries}
As is well known, the Riemann zeta function has a number of zeros, along the negative real axis at even integer points and along the critical line $x=\frac{1}{2}$. These latter appear at points whose exact positions are not known beforehand but require a lot of numerical work, the points being of a transcendental nature. Every zero in the upper half of the complex plane has a pair mirroring on the negative half plane. This work is now focusing in these nontrivial zeros. The following treatment is relying on the validity of the Riemann hypothesis \cite{Titchmarsh1999}. It requires all nontrivial zeros to reside on the critical line $x=\frac{1}{2}$ and none outside it. We are living the times of seeing the hypothesis finally proven, if not already done.

Now what will happen if a zeta function approaching a zero is divided by another zeta approaching a corresponding pairing zero on the negative half plane? As both of them go to zero one would expect the ratio to become singular. While investigating numerically the behavior of the Riemann zeta function approaching some of its zeros on the critical line, it was noted that the phase behavior of the function has a particular feature. When the argument $s$ approaches any zero $s_0$ of the zeta function $\zeta(s)$. The ratio seems to be
\begin{equation}
\lim_{s\rightarrow{s_0}}\frac{\zeta(s)}{\zeta(\bar{s})}=e^{i\theta}  \label{eqn2}
\end{equation}
with $s=\frac{1}{2}+i{t}$. One can evaluate this ratio while approaching the zero at $s_0$ from any direction and it appears that the ratio will not be singular. It can be proven as done in the following section. This leads one to think that there could exist another way to study the zeros of the zeta function. 

\section{Elimination of the Zeta Function}
Since the zeta function is real on the real axis and meromorphic elsewhere, one has according to Schwarz's reflection principle
\begin{equation}
\zeta(\bar{s})=\bar{\zeta(s)}  \label{eqn30}
\end{equation}
The function has a general form 
\begin{equation}
\zeta(s)=e^{i\phi}\rho  \label{eqn40}
\end{equation}
over the complex plane (with $\phi,\rho \in \textbf{\large{R}}$). Especially one does have along the critical line
\begin{equation}
\zeta(s_0)=e^{i\phi_0}\rho_0  \label{eqn42}
\end{equation}
and then one has the following
\begin{equation}
\bar{\zeta(s_0)}={{\rho}_0}e^{-i\phi_0}  \label{eqn52}
\end{equation}
But according to equation (\ref{eqn30}) this is equal to
\begin{equation}
\zeta(\bar{s_0})={{\rho}_0}e^{-i\phi_0}  \label{eqn56}
\end{equation}
Therefore,
\begin{equation}
\frac{\zeta(s_0)}{\zeta(\bar{s_0})}=e^{2i\phi_0}  \label{eqn60}
\end{equation}
proving the assertion and verifying the original numerical observation. This holds along the entire critical line. The function ratio (\ref{eqn2}) can be processed as follows at or very near a zero $s_0$
\begin{equation}
\frac{\zeta(s_0)}{\zeta(\bar{s_0})}=\frac{\zeta(s_0)}{\zeta(1-s_0)}  \label{eqn20}
\end{equation}
Approaching the zero can be done along the critical line from either direction. Then the ratio's amplitude remains at unity and the phase angle alone is varying. It can be done from a close distance from other directions on the complex plane while the amplitude is not unity and the phase angle is varying. The correct value is finally obtained at the zero and the amplitude reaches unity. 

The functional equation of the Riemann zeta function with the argument $s \in \textbf{\large{C}}$ is well known \cite{Riemann1858}, \cite{Siegel1932}, \cite{Edwards2001}, \cite{Titchmarsh1999}. One applies it in the following form.
\begin{equation}
\frac{\zeta(s)}{\zeta(1-s)}=\frac{(2\pi)^s}{2cos(\frac{\pi{s}}{2})\Gamma(s)}  \label{eqn71}
\end{equation}
From this it follows that at any nontrivial zero along the critical line
\begin{equation}
\frac{\zeta(s_0)}{\zeta(\bar{s_0})}=\frac{(2\pi)^{s_0}}{2cos(\frac{\pi{s_0}}{2})\Gamma(s_0)}=e^{2i\phi_0+2\pi{iN}}  \label{eqn100}
\end{equation}
This is true for any point $s_0$ on the complex plane and the right side equality is to be used in the following. One has eliminated the zeta function from this equation. One has implemented the $Mod(2\pi{i}), N\in\textbf{\large{N}}$ to be carried further on the way to the final expressions. The Weierstrass formula for the Gamma function $\Gamma(s)$ is valid for $s \in{\textbf{\large{C}}}$ 
\begin{equation}
\frac{1}{\Gamma(s)}={s}e^{\gamma{s}}\prod_{k=1}{(1+\frac{s}{k})e^{-\frac{s}{k}}}  \label{eqn120}
\end{equation}
and one can substitute it getting
\begin{equation}
e^{2i\phi_0+2\pi{iN}}=\frac{{(2\pi)}^{s_0}{s_0}e^{\gamma{s_0}}\prod_{k=1}{(1+\frac{s_0}{k})e^{-\frac{s_0}{k}}}}{2cos(\frac{\pi{s_0}}{2})}  \label{eqn140}
\end{equation}
Taking the logarithm of the equation above one will obtain
\begin{equation}
2i\phi_0+2\pi{iN}={s_0}(\gamma+ln(2\pi))+ln(\frac{s_0}{2})-ln(cos(\frac{\pi{s_0}}{2}))+\sum_{k=1}{(-\frac{s_0}{k}+ln(1+\frac{s_0}{k}))}  \label{eqn160}
\end{equation}
This equation with a complex variable $s_0$ is valid for all nontrivial zeros. 

Breaking down the $s_0$ as $s_0=\frac{1}{2}+i{t}$ where $t \in \textbf{\large{R}}$, and substituting it will lead to
\begin{equation}
2i\phi_0+2\pi{iN}=({\frac{1}{2}+i{t}})(\gamma+ln(2\pi))+ln(\frac{\frac{1}{2}+i{t}}{2})-ln(cos(\frac{\pi({\frac{1}{2}+i{t}})}{2}))  \nonumber
\end{equation}
\begin{equation}
+\sum_{k=1}{[-\frac{\frac{1}{2}+i{t}}{k}+ln(1+\frac{\frac{1}{2}+i{t}}{k})]}  \label{eqn200}
\end{equation}
The complex logarithms and the cosine term can further be broken down to real and imaginary parts as well and then one will get
\begin{equation}
2i\phi_0+2\pi{iN}=({\frac{1}{2}+i{t}})(\gamma+ln(2\pi))+ln(\frac{\sqrt{1+4{t^2}}}{4})-ln(\sqrt{\frac{cosh^2(\frac{\pi{t}}{2})+sinh^2(\frac{\pi{t}}{2})}{2}})  \nonumber
\end{equation}
\begin{equation}
+i\cdot{artan(2{t})}+i\pi{M}-i\cdot{artan(-tanh(\frac{\pi{t}}{2}))}+i\pi{L}  \nonumber
\end{equation}
\begin{equation}
+\sum_{k=1}{[-\frac{1}{2k}-\frac{it}{k}+ln(\sqrt{(1+\frac{1}{2k})^2+\frac{t^2}{k^2}})+i\cdot{artan(\frac{t}{k+\frac{1}{2}})}]}+i\pi{K}  \label{eqn220}
\end{equation}
Obviously, the real and imaginary parts of each term succeeded to get a linear form allowing a trivial separation. 
\subsection{The Real Part}
For studying the zeros, one would be interested in the real part, equation (\ref{eqn220})
\begin{equation}
0=\frac{1}{2}(\gamma+ln(2\pi))+ln(\frac{\sqrt{1+4{t^2}}}{4})+  \nonumber
\end{equation}
\begin{equation}
-ln\sqrt{\frac{cosh^2(\frac{\pi{t}}{2})+sinh^2(\frac{\pi{t}}{2})}{2}}+  \nonumber
\end{equation}
\begin{equation}
+\sum_{k=1}{[-\frac{1}{2k}+\frac{1}{2}ln(1+\frac{1}{k}+\frac{1}{4k^2}+\frac{t^2}{k^2})]}  \label{eqn250}
\end{equation}
This can be developed further simplifying to
\begin{equation}
\gamma+ln(\frac{\pi}{4})=ln(\frac{cosh(\pi{t})}{1+4{t^2}})  \nonumber
\end{equation}
\begin{equation}
-\sum_{k=1}{ln[e^{-\frac{1}{k}}({1+\frac{1}{k}+\frac{1}{4k^2}+\frac{t^2}{k^2}})]}  \label{eqn270}
\end{equation}
In spite of its apparent complexity, this equation is an identity which is valid for all $t$.
\subsection{Solving the Phases}
The imaginary part of the equation (\ref{eqn220}) can be used for getting a general expression for the limit phase $\phi_0$ at a point along the critical line
\begin{equation}
\phi_0+\frac{\pi{M}}{2}=\frac{t}{2}(\gamma+ln(2\pi))+\frac{1}{2}artan(2{t})-\frac{1}{2}{artan(-tanh(\frac{\pi{t}}{2}))}  \nonumber
\end{equation}
\begin{equation}
+\frac{1}{2}\sum_{k=1}{[\frac{-t}{k}+artan(\frac{t}{k+\frac{1}{2}})]}  \label{eqn320}
\end{equation}
The result has become $Mod({\frac{\pi}{2}}), M\in\textbf{\large{N}}$. The resulting function is odd with respect to $t$.
\subsection{Numerical Evaluation of the Phase}
\begin{figure}[htbp]
	\centering
		\includegraphics[width=1.100\textwidth]{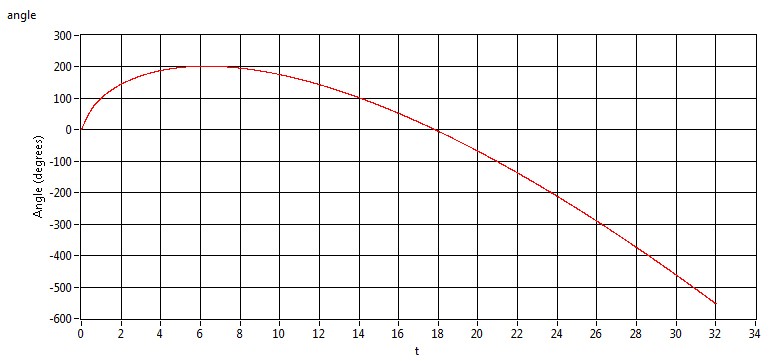}
	\caption{Phase in degrees along $x=0.5$}
	\label{Fig.1}
\end{figure}

\begin{figure}[htbp]
	\centering
		\includegraphics[width=1.100\textwidth]{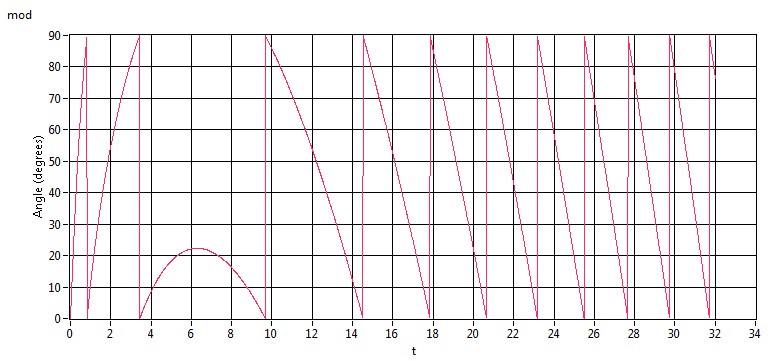}
	\caption{$Mod({\frac{\pi}{2}})$ of the phase in degrees along $x=0.5$}
	\label{Fig.2}
\end{figure}

A crude calculation of the ratio from equation (\ref{eqn60}) shows a value of -80.95 degrees $\pm{90}$ degrees for the first nontrivial zero and for the second one 77.36 degrees. One can plot the curve as a function of $t$ along the axis $x=0.5$ as in Figure  \ref{Fig.1}.
This is the angle $\phi_0$ as such. The $Mod({\frac{\pi}{2}})$ of the angle is more natural and produces Figure \ref{Fig.2}.
The negative side curve is of opposite sign as the function is odd.
\section{Discussion}
A few interesting results have been obtained. The first is that according to equation (\ref{eqn60}) the ratio between the zeta function at any zero and its conjugate is not singular, but always with unity absolute value and with a particular phase. This fact holds along the whole critical line ($x=\frac{1}{2}$) but fails immediately outside of it. The failure is not dramatic but will cause an error which increases while the point of focus is moving farther from the critical line. The second is that the equation (\ref{eqn100}) has only elementary functions left for studying the zeros of the zeta function. The third interesting finding is ensuing from the previous two, equation (\ref{eqn320}) presenting a simplified expression for calculation of phase at a zero. These equations form the main results of this work.

\end{document}